\documentclass[12pt]{amsart}
\usepackage{amsfonts,amssymb,amsthm}

\newtheorem{thm}{Theorem}
\newtheorem{cor}[thm]{Corollary}
\newtheorem{lemma}[thm]{Lemma}
\newtheorem{prop}[thm]{Proposition}

\newcommand{\R}{\mathbb{R}}
\newcommand{\C}{\mathbb{C}}
\newcommand{\E}{\mathbb{E}}
\renewcommand{\P}{\mathbb{P}}
\newcommand{\N}{\mathbb{N}}
\newcommand{\M}{\mathbb{M}}

\newcommand{\uin}{|||}
\DeclareMathOperator{\conv}{conv}
\DeclareMathOperator{\dist}{dist}
\DeclareMathOperator{\supp}{supp}
\DeclareMathOperator{\Tr}{Tr}

\title[Norms and eigenvalues of random matrices]{Concentration of norms
	and eigenvalues of random matrices}

\author[M. Meckes]{Mark W. Meckes}

\address{Department of Mathematics, Case Western Reserve University,
    Cleveland, Ohio 44106.}

\email{mwm2@po.cwru.edu}

\thanks{Partially supported by a grant from the National Science
	Foundation.}

\date{18 November 2002}

\begin{document}

\begin{abstract}
We prove concentration results for $\ell_p^n$ operator norms of
rectangular random matrices and eigenvalues of self-adjoint random
matrices. The random matrices we consider have bounded entries which are
independent, up to a possible self-adjointness constraint. Our results are
based on an isoperimetric inequality for product spaces due to Talagrand.
\end{abstract}

\maketitle


\section{Introduction} \label{intro}

In this paper we prove concentration results for norms of rectangular random
matrices acting as operators between $\ell_p^n$ spaces, and eigenvalues of
self-adjoint random matrices. Except for the self-adjointness condition when
we consider eigenvalues, the only assumptions on the distribution of the 
matrix entries are independence and boundedness.
Our approach is based on a powerful isoperimetric inequality
for product probability spaces due to Talagrand \cite{Tal1}.

Throughout this paper $X = X_{m,n}$ will stand for an $m \times n$ random
matrix with real or complex entries $x_{jk}$. (Specific technical conditions
on the $x_{jk}$'s will be introduced as needed for each result below.) 
If $1 \le p, q \le \infty$ and $A$ is an $m \times n$ matrix, we
denote by $\| A \|_{p \to q}$ the operator norm of 
$A : \ell_p^n \to \ell_q^m$. We denote by $p' = p/(p-1)$ the
conjugate exponent of $p$. For a real random variable $Y$ we denote by $\E
Y$ the expected value and by $\M Y$ any median of $Y$. Our first main
result is the following.

\begin{thm} \label{l_p-operators}
Let $1 < p \le 2 \le q < \infty$. Suppose the entries $x_{jk}$ of $X$ are
independent complex random variables, each supported in a set of diameter
at most $D$. Then
\begin{equation} \label{l_p-ineq}
\P \Bigl[ \bigl| \ \| X \|_{p \to q} - \M \| X \|_{p \to q} \bigr| \ge t
	\Bigr]
    \le 4 \exp \left[-\frac{1}{4} \left(\frac{t}{D}
    \right)^r \right]
\end{equation}
for all $t > 0$, where $r = \min \{ p', q \}$.
\end{thm}

To prove Theorem \ref{l_p-operators} we show that Talagrand's
isoperimetric inequality, which at first appears adapted primarily to prove
normal concentration for functions which are Lipschitz with respect to a
Euclidean norm, actually implies sometimes stronger concentration for
functions which are Lipschitz with respect to more general norms. In
particular, as we show in Corollary \ref{l_p-con} below, one obtains
concentration of the kind in \eqref{l_p-ineq} for convex
functions which are Lipschitz with respect to $\ell_r$ norms for $r \ge 2$.
Since such functions are automatically Lipschitz with respect to the
Euclidean norm, one can apply the known $r = 2$ case of this fact directly,
but would then obtain the upper bound with $r$ replaced by 2 in the r.h.s.~
of \eqref{l_p-ineq}. Since the conclusion of Theorem \ref{l_p-operators} is
trivial when $t / D \le 1$, the estimate \eqref{l_p-ineq} is stronger than
the estimate one would obtain this way.

To put Theorem \ref{l_p-operators} in perspective, we consider the particular
case in
which $p = q'$, $m = n$, and $\P[x_{jk} = 1] = \P[x_{jk} = -1] = 1/2$ for
all $j,k$. In this situation, 
$$n^{1/q} \le \M \| X \|_{q' \to q} \le C n^{1/q}$$
where $C > 0$ is a universal constant. Theorem \ref{l_p-operators} implies
that, while $\| X \|_{q' \to q}$ achieves values as large as $n^{2/q}$,
it is comparable to its median except on a set whose probability decays
exponentially quickly as $n \to \infty$. Furthermore, in this situation
the estimate in \eqref{l_p-ineq} is sharp as long as $n^{-1/q} t$ is 
sufficiently
large and $n^{-2/q} t$ is sufficiently small. These observations apply in
more general situations; see the remarks in Section \ref{rm} following
the proof of Theorem \ref{l_p-operators}.

\medskip

If $A$ is a self-adjoint $n \times n$ matrix, we denote by $\lambda_1(A)
\ge \lambda_2(A) \ge \cdots \ge \lambda_n(A)$ the eigenvalues of $A$,
counted with multiplicity. Our second main result is the following.

\begin{thm} \label{eigenvalues}
Suppose $m = n$ and the entries $x_{jk}$, $1 \le j \le k \le n$, of $X$
are independent complex random variables such that:
\begin{enumerate}
\item for $1 \le j \le n$, $x_{jj}$ is real and is supported on an
interval of length at most $\sqrt{2} D$; and
\item \label{r-i} for $1 \le j < k \le n$, $x_{jk}$ is either supported on
a set of diameter at most $D$; or $x_{jk} = w_{jk} ( \alpha_{jk} + i
\beta_{jk})$, where $w_{jk} \in \C$ is a constant with  $|w_{jk}| \le 1$,
and $\alpha_{jk}, \beta_{jk}$ are independent real random variables each
supported in intervals of length at most $D$;
\end{enumerate}
and that $x_{jk} = \overline{x_{jk}}$ for $k < j$. Then
\begin{equation} \label{lambda_1-ineq}
\P \bigl[ | \lambda_1(X) - \M \lambda_1(X) | \ge t \bigr]
	\le 4 e^{-t^2 / 8 D^2}
\end{equation}
for all $t > 0$, and the same holds if $\lambda_1(X)$ is replaced by
$\lambda_n(X)$. Furthermore, for each $2 \le k \le n-1$, there exists an
$M_k \in \R$ such that
\begin{equation} \label{lambda_k-ineq}
\P \bigl[ | \lambda_k(X) - M_k | \ge t \bigr] \le
    8 \exp\left[ - \left(\frac{t}{2\sqrt{2}
    (\sqrt{k} + \sqrt{k-1}) D}\right)^2\right]
    \le 8 \exp\left[ - \frac{t^2}{32 k D^2}\right]
\end{equation}
for all $t > 0$, and the same upper bound holds if $\lambda_k(X)$ is
replaced by $\lambda_{n-k+1}(X)$ and $M_k$ by $M_{n-k+1}$.
\end{thm}

Note that Theorem \ref{eigenvalues} applies in particular to the case
of real symmetric random matrices with off-diagonal entries supported in
intervals of length $D$ and diagonal entries supported in intervals of
length $\sqrt{2}D$.

The proof of Theorem \ref{eigenvalues} is also based on Talagrand's
theorem, in this case applying it only to functions which are Lipschitz
with respect to a Euclidean norm. Theorem \ref{eigenvalues} is (up to numerical
constants) a sharpening and generalization of a result of Alon,
Krivelevich, and Vu \cite{AKV}. Their proof is also based on Talagrand's
theorem, although they apply it in a very different way. For 
perspective, we note that in the particular case in which $\P[x_{jk} = 1] = 
\P[x_{jk} = -1] = 1/2$, $\M \lambda_1(X)$ is of the order $\sqrt{n}$,
while $\lambda_1(X)$ can achieve values as large as $n$. Furthermore, in
this situation the estimate in (\ref{lambda_1-ineq}) is sharp when 
$n^{-1/2} t$ is sufficiently large and $n^{-1} t$ is sufficiently small.
(See \cite{AKV} for a discussion this point when $X$ is the adjacency
matrix of the random graph $G(n, 1/2)$.)
However, the estimate in (\ref{lambda_k-ineq}) is probably not sharp in
its dependence on $k$. See the remarks in Section \ref{rm} following the
proof of Theorem \ref{eigenvalues} for details.

We note that aside from the uniform boundedness assumption, the distributions
of the independent entries of $X$ in Theorems \ref{l_p-operators} and
\ref{eigenvalues} are completely arbitrary. In particular there is no
assumption of identical distribution of independent entries, and no
assumption about the values of their means.

We emphasize that our results are of interest as bounds for large
deviations. Beginning with the work of Tracy and Widom \cite{TW1, TW2},
which has been refined and extended in \cite{Soshnikov, Ledoux3, Aubrun},
it is known that, in typical situations, the kind of result contained in
\eqref{lambda_1-ineq}, while nontrivial, is not sharp when $t$ is of smaller
order than $\sqrt{n}$. More precisely, it has been established that 
typically, one has concentration of the largest eigenvalue of the form
\begin{equation} \label{small-dev}
\P [ \lambda_1(X) - \E \lambda_1(X) \ge t ]
	\le C \exp \left[ - \max \left\{ c_1 t^2,
	 c_2 \left(n^{1/6} t\right)^{3/2} \right\} \right]
\end{equation}
in the normalization used here, where $C, c_1, c_2 > 0$ are constants.

Talagrand's theorem was first applied in the
context of random matrices by Guionnet and Zeitouni \cite{GZ}, who used it
to prove a concentration result for the spectral measure of self-adjoint
random matrices, and who also remarked that the same methods give
concentration results for other functionals of self-adjoint matrices.
For general
discussions of applications of concentration of measure phenomena to
random matrices, see the survey \cite{DS} by Davidson and Szarek and
Section 8.5 of the book \cite{Ledoux1} by Ledoux.

In Section \ref{gen-con}, we show how to obtain concentration for Lipschitz
functions on $\ell_q$ sum spaces (and more general sums of normed spaces)
from Talagrand's isoperimetric inequality. In Section \ref{rm}, we prove
Theorems \ref{l_p-operators} and \ref{eigenvalues}, and give an
infinite-dimensional version of Theorem \ref{l_p-operators} and a
version of Theorem \ref{eigenvalues} for singular values of
rectangular matrices. We also
compare the results obtained by our methods with the corresponding results
for Gaussian random matrices obtained from the Gaussian isoperimetric
inequality.


\section{General concentration results} \label{gen-con}

We first need some notation. Let $(\Omega_1, \Sigma_1, \mu_1), \ldots,
(\Omega_N, \Sigma_N, \mu_N)$ be probability spaces, $\Omega = \Omega_1
\times \cdots \times \Omega_N$, $\P = \mu_1 \otimes \cdots \otimes \mu_N$.
For $x = (x_1, \ldots, x_N) \in \Omega$, $y = (y_1, \ldots, y_N) \in
\Omega$, $h(x,y) \in \R^N$ is defined by
$$h(x,y)_j = \left\{\begin{array}{cl} 0 & \mbox{ if } x_j = y_j, \\
    1 & \mbox{ if } x_j \ne y_j. \end{array}\right. $$
For $A \subseteq \Omega$, $x \in \Omega$, $U_A(x) = \{h(x,y) : y \in A\}
\subset \R^N$. Finally, we define the convex hull distance from $x$ to $A$
by
\begin{equation} \label{f_c}
f_c(A,x) = \inf \bigl\{| z | : z \in \conv U_A(x) \bigr\},
\end{equation}
where $| \cdot |$ is the standard Euclidean norm and $\conv$ denotes
the convex hull. Talagrand's isoperimetric inequality is the following.

\begin{thm}[Talagrand \cite{Tal1}] \label{Talagrand}
Let $(\Omega, \P)$ be a product probability space as above. For any 
$A \subseteq \Omega$,
$$\int_\Omega \exp \left[\frac{1}{4} f_c^2 (A, x) \right] d\P(x)
    \le \frac{1}{\P(A)},$$
which by Chebyshev's inequality implies
$$ \P \bigl( \{ x : f_c (A, x) \ge t \} \bigr) \le
	\frac{1}{\P (A) } e^{- t^2 / 4}$$
for all $t > 0$.
\end{thm}

As in \cite{Tal1}, we have ignored measurability issues in the
statement of Theorem \ref{Talagrand}. To be strictly correct,
the integrals and probabilities which appear must be replaced by upper
integrals and outer probabilities; however, this issue is 
irrelevant in applications, since one typically applies such a result to
estimate expressions in which all the functions and sets which appear are
measurable.

Let $\| \cdot \|_E$ be a 1-unconditional norm on $\R^N$, by which we mean
that the standard basis of $\R^n$ is a 1-unconditional basis for 
$\| \cdot \|_E$ (see \cite{LT-Banach}; such a norm is also sometimes called
absolute).
For normed vector spaces $(V_j, \| \cdot \|_{V_j})$, $j = 1, \ldots, N$,
we denote by
$$V_E = \left(\bigoplus_{j=1}^N V_j \right)_E$$
the direct sum of vector spaces with the norm
$$\| (v_1, \ldots, v_N) \|_{V_E}
	= \bigl\| \bigl( \| v_1 \|_{V_1}, \ldots, \| v_N \|_{V_N} \bigr) 
	\bigr\|_E.$$

Theorems \ref{l_p-operators} and \ref{eigenvalues} will be proved using
the following consequence of Theorem \ref{Talagrand}.

\begin{cor} \label{l_p-con}
Let $V$ be the $\ell_q$ sum of the normed vector spaces $(V_j, \|
\cdot \|_{V_j})$, $j = 1, \ldots, N$ for $q \ge 2$;
that is, $V = V_{\ell_q^N}$ in
the notation above. For $j = 1, \ldots, N$, let $\mu_j$ be a probability
measure on $V_j$ which is supported on a compact set of diameter at most 1.
Let $\P = \mu_1 \otimes \cdots \otimes \mu_n$. Suppose $F : V \to \R$ is
1-Lipschitz and quasiconvex, that is, $F^{-1} \bigl((-\infty,a] \bigr)$
is convex for all $a \in \R$. Then
\begin{equation} \label{l_p-sum-ineq}
\P\bigl[| F - \M F | \ge t \bigr] \le 4 e^{- t^q / 4}
\end{equation}
for all $t > 0$.\end{cor}

We will postpone the proof of Corollary \ref{l_p-con} until after some
remarks.
The $q = 2$ case of Corollary \ref{l_p-con} has been widely noted and 
applied in various degrees of generality already; see \cite{Tal1, Tal2,
Ledoux1, Schecht} and their references. 
As observed in the introduction, if $F$ is 1-Lipschitz
with respect to the $\ell_q$ sum norm on $V$, then $F$ is also 1-Lipschitz
with respect to the $\ell_2$ sum norm on $V$. However, applying the $q = 2$
case of Corollary \ref{l_p-con} directly in this situation yields only the
weaker upper bound $4 e^{- t^2 / 4}$ in the inequality \eqref{l_p-sum-ineq}.

Corollary \ref{l_p-con} can be applied in the
case that $F$ is $L$-Lipschitz and each $\mu_i$ is supported on a set of
diameter at most $D$, by replacing $t$ with $t / LD$ in the r.h.s.~ of
\eqref{l_p-sum-ineq}. This fact is used
implicitly in the proofs in Section \ref{rm}.
The conclusion of Corollary \ref{l_p-con} also holds when $F$ is
replaced by $-F$, so that Proposition \ref{med-con} also applies when $F$ is
quasiconcave, that is, when $F^{-1} \bigl([a,\infty)\bigr)$
is convex for all $a \in
\R$. In particular, Corollary \ref{l_p-con} applies to both convex and
concave Lipschitz functions.
Talagrand gives an example which shows that some form of convexity
assumption in Corollary \ref{l_p-con} cannot be removed in general in
\cite{Tal3}, in which the special case of Theorem \ref{Talagrand} for the
uniform measure on the discrete cube was first proved.

The conclusion of Corollary \ref{l_p-con} does not hold in general for
functions which are only Lipschitz with respect to an $\ell_p$ sum norm
for $1 \le p < 2$
without the introduction of dimension dependent constants, even if the
bound $4e^{-t^p / 4}$ is replaced by any other
dimension independent function which approaches 0 at infinity.
To see this, let $\{ Y_j : j \in \N \}$ be independent random
variables with $\P[Y_j = 0] = \P[Y_j = 1] = 1/2$ for each $j$, and $S_n =
\sum_{j=1}^n Y_j$ for an arbitrary $1 \le n \le N$.
Then $S_n^{1/p} = \|(Y_1, \ldots, Y_n)\|_p$, and
$(n/2)^{1/p}$ is a median for $S_n^{1/p}$. Suppose we have a concentration
result which implies there exists a function $f$ with $\lim_{t \to \infty}
f(t) = 0$ such that for all $n$ and for all $t > 0$,
\begin{equation}
\P \bigl[ S_n^{1/p} - \M S_n^{1/p} \ge t \bigr] \le f(t). \label{f(t)}
\end{equation}
Then by Taylor's theorem applied at $t = 0$,
\begin{equation*} \begin{split}
\P\left[\frac{S_n-n/2}{\sqrt{n}/2} \ge t\right] & =
    \P\left[S_n^{1/p} \ge \left(\frac{n}{2} +
    \frac{\sqrt{n}}{2} t \right)^{1/p} \right] \\
& = \P\left[S_n^{1/p} - \left(\frac{n}{2}\right)^{1/p} \ge
    \frac{1}{2^{1/p}p} n^{\frac{1}{p} - \frac{1}{2}}t +
    O\left(n^{\frac{1}{p}-1}\right)t^2\right] \\
& \le f\left(\frac{1}{2^{1/p}p}n^{\frac{1}{p}-\frac{1}{2}}t
    +O\left(n^{\frac{1}{p}-1}\right)t^2\right),
\end{split} \end{equation*}
which implies that for any $t > 0$,
$$\lim_{n \to \infty}
    \P\left[\frac{S_n - n/2}{\sqrt{n}/2} \ge t \right] = 0,$$
which contradicts the central limit theorem.  Since concentration of the
kind in the inequality (\ref{f(t)}) about any value implies concentration
about a median (with a possibly different function $f$), no such
concentration result holds for $S_n^{1/p}$ when $1 \le p < 2$.

It is also not difficult to see that for the examples of $S_n^{1/q}$ with
$q \ge 2$, the concentration result of Corollary \ref{l_p-con} is sharp,
up to the values of numerical constants, when $c = c(q) \le n^{-1/q} t
\le 1 - 2^{-1/q}$. Moreover, by gluing together copies of $S_n^{1/q}$ for
different values of $n$, one obtains an example of a Lipschitz function
for which
Corollary \ref{l_p-con} is sharp for the entire nontrivial range of $t$.

It is more typical to state results of the type
in Corollary \ref{l_p-con} in terms of deviations of a random variable
from the mean rather than the median. This difference is inessential,
since this level of concentration implies that the median and mean cannot
be too far apart. For example, in the situation of Corollary \ref{l_p-con},
we have
\begin{equation*} \begin{split}
| \E F - \M F| & \le \E |F - \M F|
    = \int_0^\infty \P \bigl[ |F - \M F| \ge t \bigr] dt \\
& \le 4 \int_0^\infty e^{-t^q / 4} dt 
	= 4^{1+ \frac{1}{q}} \Gamma \bigl(1 + \tfrac{1}{q} \bigr).
\end{split} \end{equation*}

We now turn to the proof of Corollary \ref{l_p-con}. Rather than prove
Corollary \ref{l_p-con} directly from Theorem
\ref{Talagrand}, we will deduce it as a special case of Proposition
\ref{med-con} below, which uses Theorem \ref{Talagrand} to derive
concentration for functions which are Lipschitz with respect to an arbitrary
1-unconditional norm, in terms of a kind of modulus function for the norm.
Theorem \ref{Talagrand} bounds the size of the set of points which are far
from a set $A$ in terms of the convex hull distance $f_c (A, \cdot )$ 
from $A$. Thus it
provides concentration for functions which satisfy a Lipschitz type
condition with respect to the convex hull distance. However, since in general
this distance is not induced by a metric, some care is needed in its
application.

Let $\| \cdot \|_E$ be a 1-unconditional norm on $\R^N$ as above. We define
$$K_E (t) = \inf \bigl\{ | x | : \| x \|_E \ge t, \, \| x \|_\infty \le 1 
	\bigr\},$$
where we use the convention that $\inf \emptyset = \infty$.

\begin{prop} \label{med-con}
Let $V_E$ be as described before Corollary \ref{l_p-con} and let $K_E$ be as
above. For $j = 1, \ldots, N$, let $\mu_j$ be a probability
measure on $V_j$ which is supported on a compact set of diameter at most 1.
Let $\P = \mu_1 \otimes \cdots \otimes \mu_n$. Suppose $F : V_E \to \R$ is
1-Lipschitz and that $F$ is quasiconvex. Then
\begin{equation} \label{con-ineq}
\P \bigl[ | F - \M F | \ge t \bigr] \le 4 \exp \left[ - \frac{1}{4} 
	\bigl( K_E(t) \bigr)^2 \right]
\end{equation}
for all $t > 0$.
\end{prop}

It is easy to verify that $K_{\ell_q^N}(t) \ge t^{q/2}$ for $q \ge 2$ and
any $N \in \N$ (c.f.
Lemma \ref{K_E} below), so that Corollary \ref{l_p-con} follows immediately
from Proposition \ref{med-con}.

\begin{proof}[Proof of Proposition \ref{med-con}]
First we show that
\begin{equation} \label{dist_E}
K_E \bigl(\dist (x, \conv A)\bigr) \le f_c (A, x)
\end{equation}
for $x = (x_1, \ldots, x_N) \in \supp (\P)$ and $\emptyset \ne A \subseteq
V_E$, where $\dist$ is the distance in the normed space $V_E$.
Let $y^k = (y^k_1, \ldots, y^k_N) \in A \cap \supp (\P)$ and $0 \le
\theta_k \le 1$ for $k = 1, \ldots, n$ such that $\sum_{k=1}^n \theta_k = 1$.
Then for each $j = 1, \ldots, N$,
\begin{equation*}
\left\| x_j - \sum_{k=1}^n \theta_k y^k_j \right\|_{V_j}
	\le \sum_{k=1}^n \theta_k \bigl\| (x_j - y^k_j) \bigr\|_{V_j}
	\le \sum_{k=1}^n \theta_k h(x, y^k)_j
\end{equation*}
since $x_j, y_j^k \in \supp (\mu_j)$ for each $j, k$. Then by 
unconditionality,
$$\dist (x, \conv A) \le
	\left\| x - \sum_{k=1}^n \theta_k y^k \right\|_{V_E}
	\le \left\| \sum_{k=1}^n \theta_k h(x, y^k) \right\|_E,$$
and so
\begin{equation*}
K_E \bigl(\dist (x, \conv A)\bigr) \le \left| \sum_{k=1}^n \theta_k h(x, y^k)
	\right|.
\end{equation*}
The inequality (\ref{dist_E}) now follows since the r.h.s.~ of (\ref{f_c}) is
precisely the infimum of this last expression over all such finite
sequences $y^k, \theta_k$, $k = 1, \ldots, n$. Therefore, by Theorem
\ref{Talagrand}, for any $A \subseteq V_E$,
$$ \P (A) \P \bigl( \bigl\{ x : K_E \bigl(\dist (x, \conv A) \bigr) \ge t 
	\bigr\} \bigr) \le e^{- t^2 / 4}$$
for all $t > 0$. Thus if $F$ is quasiconvex and 1-Lipschitz on $V_E$, we
have that for any $a \in \R$, $t > 0$,
\begin{equation*} \begin{split}
\P[F \le a] \P[F \ge a + t]
& \le \P[F \le a] 
	\P \left(\left\{ x : K_E \bigl( \dist(x , 
	F^{-1} \bigl((-\infty, a]\bigr)\bigl) \ge K_E(t) \right\}\right) \\
& \le \exp \left[ -\frac{1}{4} \bigl(K_E(t)\bigr)^2 \right].
\end{split} \end{equation*}
Applying this in turn with $a = \M F$ and $a = \M F - t$, we get
$$\P[F - \M F \ge t] \le 2 \exp \left[ - \frac{1}{4} 
	\bigl(K_E(t)\bigr)^2 \right],$$
$$\P[F - \M F \le -t] \le 2 \exp \left[ - \frac{1}{4} 
	\bigl(K_E(t)\bigr)^2 \right],$$
for every $t > 0$.
\end{proof}

In order to apply Proposition \ref{med-con}, one needs to estimate
the function $K_E$. This is of most interest if one can bound $K_{E_j}$
uniformly for some family of spaces $E_j$ for which $\sup_j \dim (E_j)
= \infty$. This is not difficult to do for certain classes of
spaces. For a nonincreasing sequence $w = (w_1, w_2, \ldots )$ of positive
numbers and $p \ge 1$, the $N$-dimensional Lorentz space $\ell_{w, p}^N$ is
$\R^N$ with the norm
\begin{equation*} 
\| x \|_{w, p} = \left( \sum_{j=1}^N w_j a_j^p \right)^{1/p},
\end{equation*}
where $\{ a_j : 1 \le j \le N \}$ is the nonincreasing rearrangement of
$\{ | x_j | : 1 \le j \le N \}$. For an Orlicz function $\psi$, that is,
a convex nondecreasing function $\psi : \R_+ \to \R_+$ such that 
$\psi(0) = 0$ and $\lim_{t \to \infty}\psi(t) = \infty$, the $N$-dimensional
Orlicz space $\ell_\psi^N$ is $\R^N$ with the norm
\begin{equation*} 
\| x \|_\psi = \inf \left\{ \rho > 0 : \sum_{j=1}^N 
	\psi \left( \frac{|x_j|}{\rho} \right) \le 1 \right\}.
\end{equation*}
Observe that $\ell_p^N = \ell_{w, p}^N$ if $w_j = 1$ for $j = 1, \ldots, N$, and
$\ell_p^N = \ell_\psi^N$ if $\psi(t) = t^p$, $p \ge 1$. For these two
classes of spaces, we have the following elementary estimates, which we
state without proof.

\begin{lemma} \label{K_E}
If $p \ge 1$ and $w \in \ell_r$ for some $r$ such that 
$\max \{1, 2/p \} \le r' < \infty$, then
$$ K_{\ell_{w, p}^N} (t) \ge \| w  \|_r^{-r' / 2} t^{ p r' / 2}.$$
If $\psi$ is any Orlicz function, then
$$ K_{\ell_\psi^N} (t) \ge \inf_{0 < u \le 1} \frac{u}{\sqrt{\psi(u / t)}}.$$
In particular, $K_{\ell_q^N}(t) \ge t^{q / 2}$ for $q \ge 2$.
\end{lemma}

Note that the estimates in Lemma \ref{K_E} may be trivial and are not
necessarily optimal, but when they are nontrivial, they are valid in
all dimensions. By
considering vectors $x \in \{ 0, 1 \}^N$, one can see that the estimate
$K_{\ell_q^N}(t) \ge t^{q/2}$ for $q \ge 2$ is sharp for $t = k^{1/q}$,
$k = 1, \ldots, N$.

Observe that the proof of Proposition \ref{med-con} actually
gives separate tail estimates for deviations of $F$ above and below its
median; the same is therefore true of Corollary \ref{l_p-con} as well.
The full generality of Proposition \ref{med-con} can in fact be
derived with some amount of argument from the (known) $q = 2$ case of
Corollary 
\ref{l_p-con}, using these bounds separately; however, we find it
simpler to argue directly from the isoperimetric inequality of Theorem
\ref{Talagrand} as above. One could alternatively prove Corollary
\ref{l_p-con} by proving an $\ell_q$ version of Theorem \ref{Talagrand},
by defining an $\ell_q$ convex hull distance $f_q(A, x) = \inf \{ 
\| z \|_q : z \in \conv U_A(x) \}$ and mimicking the proof of Theorem
\ref{Talagrand}; or as a corollary to the more general and abstract
Theorem 4.2.4
in \cite{Tal1}. However, this approach would result only in a slight 
sharpening of the constant $1/4$ which appears in the exponent.

We remark that to use Proposition \ref{med-con} to full advantage for
non-Euclidean norms, one must use a nonlinear lower bound on $K_E$ and
make use of the restriction $\| x \|_\infty \le 1$. If, for example, one
uses only the fact that $\| x \|_q \le | x |$ for all $x$ when $q \ge 2$,
then one is using no more than the fact that a function which is 
1-Lipschitz with respect to the $\ell_q$ norm is 1-Lipschitz with respect
to the $\ell_2$ norm, which, as we have observed already in the introduction,
leads to a weaker concentration result.

\medskip

It is instructive to compare the general concentration results above and
the applications in the next section with the corresponding results for
Gaussian measures. We begin by recalling the functional form of the
Gaussian isoperimetric inequality, due independently to Borell
\cite{Borell} and Sudakov and Tsirel'son \cite{ST}. Let $\gamma_N$ be the
standard Gaussian measure on $\R^N$ defined by $d \gamma_N(x) = (2
\pi)^{-N / 2} e^{-|x|^2 / 2} dx$, where $| \cdot |$ is again the standard
Euclidean norm.

\begin{thm}[Borell, Sudakov-Tsirel'son] \label{Gaussian}
Let $F : \R^N \to \R$ be 1-Lipschitz with respect to the Euclidean metric
on $\R^N$. Then
$$\gamma_N \bigl(\{x : F(x) \ge \M F + t \} \bigr)
    \le 1 - \gamma_1 \bigl((-\infty, t] \bigr) < \frac{1}{2} e^{-t^2 / 2}$$
for all $t > 0$.
\end{thm}

Observe that by composing $F$ with an affine contraction, one obtains the
same conclusion in Theorem \ref{Gaussian} if the standard Gaussian measure
$\gamma_N$ is replaced by the product of one-dimensional Gaussian measures
with arbitrary means and variances at most 1. Thus the $q = 2$ case of
Corollary \ref{l_p-con} provides a level of concentration for quasiconvex
Lipschitz functions of independent bounded random variables comparable to
the concentration of Lipschitz functions of independent Gaussian random
variables with bounded variances.

A similar concentration principle is obeyed by any
probability measure which satisfies a logarithmic Sobolev inequality (see
\cite{Ledoux2}). Specifically, if $\mu$ is a probability measure on
$\R^N$ which has logarithmic Sobolev constant at most 1, and $F : \R^N \to
\R$ is 1-Lipschitz with respect to the Euclidean metric on $\R^N$, then
$$\mu \bigl(\{ x : F(x) \ge \E F + t \} \bigr) \le e^{-t^2 / 2}$$
for all $t > 0$. Since logarithmic Sobolev inequalities tensorize, one
obtains concentration for Lipschitz functions of independent random
variables whose distributions have uniformly bounded logarithmic Sobolev
constants. In particular, whenever we state concentration results below
for random matrices with Gaussian entries, similar results hold under the
weaker assumption of entries with uniformly bounded logarithmic Sobolev
constants. We remark that Guionnet and Zeitouni \cite{GZ} also proved a
concentration result for the spectral measure in the case that the matrix
entries satisfy a logarithmic Sobolev inequality.


\section{Norms and eigenvalues of random matrices} \label{rm}

Since any norm on a real or complex vector space is a convex function,
Proposition \ref{med-con} can be applied directly to obtain concentration of
norms of a random matrix $X$; all that is necessary is to estimate the
function $K_E$, or the Lipschitz constant of the given norm with respect
to one for which a bound on $K_E$ is known. Note that by the triangle
inequality,
$$ \bigl| \ \| x \| - \| y \| \ \bigr| \le \| x - y \| $$
for any norm, which implies that to estimate the Lipschitz constant of one
norm with respect to another norm, it suffices to estimate the appropriate
equivalence constant.

\begin{proof}[Proof of Theorem \ref{l_p-operators}]
For an $m \times n$ matrix $A$, let $A_j \in \C^n$ denote the $j^{\rm th}$
row of $A$. Then H\"older's inequality implies
\begin{equation} \label{l_p-op-ineq}
\| A \|_{p \to q} \le \bigl\| \bigl( \| A_1 \|_{p'}, \ldots, \| A_m \|_{p'}
	\bigr) \bigr\|_{q} \le \| (a_{jk}) \|_r,
\end{equation}
where $(a_{jk})$ represents the matrix $A$ thought of as an element of
$\C^{mn}$, and we recall that $r = \min \{p', q\}$. The claim follows
by using this estimate and taking $V_j = \C$ for each $j$ in Corollary
\ref{l_p-con}. Alternatively, the inequality \eqref{l_p-op-ineq} and
Lemma \ref{K_E} imply that
$$K_{L(\ell_p^n, \ell_q^m)} (t) \le t^{r/2},$$
where $L(\ell_p^n, \ell_q^m)$ is identified with $\C^{mn}$ via the
standard bases, so that the claim follows from Proposition \ref{med-con}.
\end{proof}

We remark that Theorem \ref{l_p-operators} can be extended to more general
norms on ${\mathfrak M}_{m, n}( \C )$ by using Proposition \ref{med-con}
together with estimates on the corresponding function $K_E$. In particular,
as long as one has the appropriate Lipschitz estimates, the underlying normed
spaces need not be unconditional, nor must the norm on matrices even be
an operator norm.

Now for comparison, we let $G = G_{mn}$ be an $m \times n$ random matrix whose
entries are
independent Gaussian random matrices with arbitrary means and variances at
most 1. For $1 \le p \le 2 \le q \le \infty$, $\| A \|_{p \to q} \le \| A
\|_2$ for any $m \times n$ matrix $A$, where $\| A \|_2$ is the
Hilbert-Schmidt norm of $A$. Then Theorem \ref{Gaussian} implies that
$$\P \Bigl[ \bigl| \ \| G \|_{p \to q} - \M \| G \|_{p \to q} \bigl| \ge t 
	\Bigl] < e^{-t^2 / 2}$$
for all $t > 0$. Observe that this is comparable to what one would obtain
in the cases of independent bounded entries by using only the $q = 2$ case
of Corollary \ref{l_p-con}.

Theorem \ref{l_p-operators} implies that the order of fluctuations
of $\| X_{m,n} \|$ about its median is $O(1)$, independent of $m$ and $n$. In
typical situations, the median itself grows without bound as $m$ or $n$
does. Suppose for example that $\E |x_{jk}| \ge c > 0$ for all $j, k$. (In
the situation of Theorem \ref{l_p-operators}, this will be the case if
each $x_{jk}$ is real, $|x_{jk}| \le 1$, $\E x_{jk} = 0$, and $x_{jk}$
has variance at least $c$.) Then
\begin{equation*} \begin{split}
\E \| X \|_{p \to q} & \ge \E \| X e_1 \|_q
    \ge m^{1/q - 1} \E \| X e_1 \|_1 \\
& = m^{1/q - 1} \sum_{j=1}^m \E |x_{j1}| \ge c m^{1/q}.
\end{split} \end{equation*}
Since $\| X \|_{p \to q} = \| X^* \|_{q' \to p'}$, we obtain $\E \| X
\|_{p \to q} \ge c \max \{ m^{1/q}, n^{1/p'} \}$. As remarked earlier,
$\M \| X \|_{p \to q}$ will also have at least this order when the 
hypotheses of Theorem \ref{l_p-operators} are satisfied.

A similar upper estimate is possible in the case $p = q'$. Suppose that
each $x_{jk}$ is a symmetric real random variable such that $|x_{jk}| \le
1$. We note first that by the Riesz convexity theorem,
$$ \| X \|_{q' \to q} \le \| X \|_{2 \to 2}^{\frac{2}{q}}
	\| X \|_{1 \to \infty}^{1 - \frac{2}{q}}
	\le \| X \|_{2 \to 2}^{\frac{2}{q}}.$$
By the contraction principle (see \cite[Theorem 4.4]{LT}),
$$ \E \| X \|_{2 \to 2} \le \E \| \tilde{X} \|_{2 \to 2},$$
where $\tilde{X} = \tilde{X}_{m,n}$ is an $m \times n$ matrix whose
entries are independent Rademacher (Bernoulli) random variables; that is,
$\P[ \tilde{x}_{jk} = 1] = \P[ \tilde{x}_{jk} = -1] = 1/2$ for all $j, k$.
By standard comparisons between Rademacher and Gaussian averages and
Chevet's inequality \cite{Chevet} (see also \cite{LT}),
$$ \E \|  \tilde{X} \|_{2 \to 2} \le C \bigl(m^{1/2} + n^{1/2} \bigr),$$
where $C > 0$ is an absolute numerical constant. Therefore in this situation,
$$ \E \| X \|_{q' \to q} \le 2 C \max \bigl\{ m^{1/q}, n^{1/q} \bigr\}.$$
(The argument above is entirely standard and the estimate is probably
known, although we could not find a reference in the literature.)

The example of $\tilde{X}$ above can be used to show that the estimate in 
Theorem
\ref{l_p-operators} is sharp for large enough values of $t$ up to numerical
constants in the case that $p = q'$. For $1 \le a \le m$, $1 \le b \le n$,
\begin{equation*}
\P \bigl[ \| \tilde{X} \|_{q' \to q} \ge (ab)^{1/q} \bigr] \ge
    	\P \bigl[ \tilde{X} \text{ has an $a \times b$ all-1 submatrix} \bigr]
	\ge 2^{-ab},
\end{equation*}
so that 
$$\P \bigl[ \| \tilde{X} \|_{q' \to q} \ge t \bigr] \ge 2^{-t^q}$$
for $t = (ab)^{1/q}$, $a = 1, \ldots, m$, $b = 1, \ldots, n$. Together with
the above upper bound on $\E \| \tilde{X} \|_{q' \to q}$, this implies that
in this situation, the concentration result of Theorem \ref{l_p-operators}
is sharp when $(\max \{m, n \})^{-1/q} t$ is sufficiently large, up to the
values of numerical constants.

For $p, q$ in other ranges, one can derive concentration for $\| X
\|_{p \to q}$ by comparing the $\ell_{p'}^n$ or $\ell_q^m$ norm to the
$\ell_2$ norm of the appropriate dimension. In this case one will obtain
concentration on a scale which depends on $m$ or $n$. For example, in the
situation of Theorem \ref{l_p-operators} one has
$$\P \Bigl[ \bigl| \ \|X\|_{p \to q} - \M \| X \|_{p \to q} \bigr| 
	\ge t \Bigr] \le 4 \exp\left[ - \frac{t^2}{4 m^{\frac{2}{q} - 1} 
	n^{\frac{2}{p'} - 1}} \right]$$
if $1 < q \le 2 \le p < \infty$.

\medskip

Since the conclusion of Theorem \ref{l_p-operators} is independent of
dimension, one can derive the following infinite dimensional version
for kernel operators from $\ell_p$ to $\ell_q$.

\begin{cor} \label{l_p-cor}
Let $1 < p \le 2 \le q < \infty$, and let $c_{jk} \ge 0$, $j, k \in \N$,
be constants such that 
\begin{equation} \label{p-summing}
\left( \sum_{j=1}^\infty \left( \sum_{k=1}^\infty c_{jk}^{p'}
	\right)^{q/p'} \right)^{1/q} < \infty.
\end{equation}
Suppose that $x_{jk}$, $j, k \in \N$
are independent complex random variables each supported in a set of
diameter at most D, such that $|x_{jk}| \le c_{jk}$ for all
$j, k$ . Define the random operator $X : \ell_p \to \ell_q$ by setting
$X(e_j) = \sum_{k=1}^\infty x_{jk} e_k$. Then
$$\P \Bigl[ \bigl| \ \|X\| - \M \| X \| \ \bigr| \ge t\Bigr] \le 
	4 \exp \left[-\frac{1}{4} \left(\frac{t}{D}
	\right)^r \right]$$
for all $t > 0$, where $\| X \|$ is the operator norm of $X$ and $r = \min
\{ p', q \}$.
\end{cor}

We remark that when $p = q'$, the l.h.s.~ of \eqref{p-summing} was shown by
Persson \cite{Persson} to coincide
with both the $q$-summing norm $\pi_q(T)$ and the $q$-nuclear norm
$\nu_q(T)$ of the kernel operator $T:\ell_{q'} \to \ell_q$ given by
$T(e_j) = \sum_{k=1}^\infty c_{jk} e_k$.

\begin{proof}[Proof of Corollary \ref{l_p-cor}]
The fact that $|x_{jk}| \le c_{jk}$ implies that $\| X \| < \infty$ always.
Apply Theorem \ref{l_p-operators} to the $n \times n$ upper-left corner
of the infinite matrix $(x_{jk})$, and use \eqref{p-summing} and the estimate
$|x_{jk}| \le c_{jk}$ to pass to the limit $n \to \infty$.
\end{proof}

Note that by taking $c_{jk} = 0$ when $j > m$ or $k > n$ in Corollary
\ref{l_p-cor}, we recover Theorem \ref{l_p-operators}, so that these
two statements are formally equivalent.

\medskip

We now specialize to the case in which $m = n$ and consider $X$ as an
operator on $\ell_2^n$, so that we use only the $q = 2$ case of
Corollary \ref{l_p-con}. Guionnet and Zeitouni \cite{GZ} were the first
to note that this concentration
theorem implies normal concentration for any function on matrices (or
self-adjoint matrices) which is convex and Lipschitz with respect to the
Hilbert-Schmidt norm. For example, we have the following. Let the entries
$x_{jk}$ of $X$ all be independent, and satisfying the condition
(\ref{r-i}) in the statement of Theorem \ref{eigenvalues}, and for 
simplicity let $D = 1$. For $1 \le p \le \infty$, we denote by 
$\| A \|_p$ the Schatten $p$-norm of an $n \times n$ matrix $A$ (see,
e.g., \cite{Bhatia}). Then for all $t > 0$,
$$\P \Bigl[ \bigl| \ \| X \|_p - \M \| X \|_p \bigr| \ge t\Bigr]
	 \le 4 e^{-t^2/4}$$
for $2 \le p \le \infty$, and
$$\P \Bigl[ \bigl| \ \| X \|_p - \M \| X \|_p \bigr| \ge t \Bigr] \le
    4 \exp \left[ - \frac{t^2}{4 n^{\frac{2}{p} - 1}}\right]$$
for $1 \le p < 2$. (In particular, we observe that when $p = q = 2$, the
conclusion of Theorem \ref{l_p-operators} holds when the matrix entries
$x_{jk}$ satisfy condition (\ref{r-i}) in the statement of Theorem
\ref{eigenvalues}.)
Furthermore, since $\uin A \uin \le \| A \|_1 \le
\sqrt{n} \| A \|_2$ for any unitarily invariant norm $\uin \cdot \uin$ on
${\mathfrak M}_n(\C)$ satisfying $\uin E_{11} \uin = 1$, it follows that
$$\P \Bigl[ \bigl| \ \uin X \uin - \M \uin X \uin \ \bigr| \ge t \Bigr]
	\le 4 e^{-t^2 / 4n}$$
for all $t > 0$ for any such norm. Each of these observations is in fact
a special case of the tail inequalities for norms of sums of independent
vector-valued random variables which were the original motivation for 
Talagrand's
development of Theorem \ref{Talagrand} and related concentration theorems.

We now consider eigenvalues of a self-adjoint random matrix. Although
these are not (except in the extreme cases) quasiconvex or quasiconcave
functions, Corollary \ref{l_p-con} can still be used to derive
concentration.

\begin{proof}[Proof of Theorem \ref{eigenvalues}]
For simplicity, we assume $D = 1$. First observe that
$$ \| X \|_2 = \left(\sum_{j,k = 1}^n |x_{jk}|^2 \right)^{1/2}
    = \sqrt{2} \left( \sum_{j=1}^n \left| \frac{x_{jj}}{\sqrt{2}}
    \right|^2
    + \sum_{1 \le j < k \le n} |x_{jk}|^2 \right)^{1/2}.$$
We suppose for simplicity that each of the upper-diagonal entries $x_{jk}$
for $j < k$ is supported in a set of diameter at most 1. (The argument is
similar in the case that for some $j < k$, $x_{jk} = w_{jk} (\alpha_{jk} +
i \beta_{jk})$ as in the statement of the theorem.) Note that
$\frac{x_{jj}}{\sqrt{2}}$, $j = 1, \ldots, n$, and $x_{jk}$, $1 \le j < k
\le n$, are independent random variables in $\R$ or $\C$, each supported
in a set of diameter at most 1. $\| X \|_2$ is $\sqrt{2}$ times the
$\ell_2$ sum norm of the direct sum of $n$ copies of $\R$ and $\binom{n}{2}$
copies of $\C$ spanned by these variables.

Recall also that
$$ \| X \|_2 = \left( \sum_{k=1}^n \lambda_k(X)^2 \right)^{1/2},$$
which implies that each $\lambda_k(X)$ is a 1-Lipschitz function of $X$
with respect to $\| X \|_2$. The first claim now follows directly from
Corollary \ref{l_p-con} with $V_j = \C$ or $V_j = \R$ for each $j$, since
$\lambda_1$ is a convex function, and $\lambda_n$ is concave.

To prove the second claim, we introduce the following functions for a
self-adjoint matrix $A$. For $k = 1, \ldots, n$, let
$$F_k(A) = \sum_{j=1}^k \lambda_j(A),$$
$$G_k(A) = \sum_{j=1}^k \lambda_{n-j+1}(A) = \Tr A - F_k(A).$$
Then $F_k$ is positively homogeneous (of degree 1), and $F_k(-A) =
-G_k(A)$. From this it follows that
$$|F_k(A) - F_k(B)| \le \max \{F_k(A-B), -G_k(A-B)\}
    \le \sqrt{k}\| A - B \|_2,$$
$$|G_k(A)-G_k(B)| \le \sqrt{k} \| A - B \|_2.$$
Moreover, $F_k$ is convex and $G_k$ is concave for each $k$; this follows
from Ky Fan's maximum principle (see, e.g., \cite{Bhatia}) or Davis's
characterization \cite{Davis} of all convex unitarily invariant functions
of a self-adjoint matrix. Let $M_k = \M F_k - \M F_{k-1}$. Then by
Corollary \ref{l_p-con}, for any $0 \le \theta \le 1$,
\begin{equation*} \begin{split}
\P \bigl[ |\lambda_k(X) - M_k | \ge t \bigr] & =
    \P \bigl[ | (F_k(X) - \M F_k(X)) - (F_{k-1}(X) - \M F_{k-1}(X)) |
	 \ge t \bigr] \\
& \le \P \bigl[ |F_k(X) - \M F_k(X) | \ge \theta t \bigr] \\
& \quad + \P \bigl[ |F_{k-1}(X) - \M F_{k-1}(X) | \ge (1 - \theta)t \bigr] \\
& \le 4 \exp\left[ -\left(\frac{\theta t}{2\sqrt{2k}}\right)^2 \right]
    + 4 \exp\left[ -\left(\frac{(1 - \theta) t}{2\sqrt{2(k-1)}}
    \right)^2 \right].
\end{split} \end{equation*}
The estimate \eqref{lambda_k-ineq} now follows by letting $\theta =
\sqrt{k} / (\sqrt{k} + \sqrt{k-1})$. (This is not the optimal value of
$\theta$, but optimizing at this point would only result in a slight 
sharpening of the constants, and not of the dependence on $t$ or $k$.)
The claim for $\lambda_{n-k+1}(X)$
follows similarly, using $G_k(X)$ in place of $F_k(X)$, or as a formal
consequence by replacing $X$ with $-X$.
\end{proof}

Now, for comparison, we let $H_n$ be an $n \times n$ random matrix with entries
$h_{jk}$, $1 \le j, k \le n$, such that:
\begin{enumerate}
\item the entries $h_{jk}$, $1 \le j \le k \le n$ are independent Gaussian
random variables,
\item the variance of $h_{jk}$ for $1 \le j < k \le n$ is at most 1,
\item the variance of $h_{jj}$ is at most $\sqrt{2}$ for $1 \le j \le n$, and
\item $h_{jk} = h_{kj}$ for $k < j$.
\end{enumerate}
Then for each $1 \le k \le n$, Theorem \ref{Gaussian} implies that
$$\P \bigl[ | \lambda_k(H_n) - \M \lambda_k(H_n) | \ge t \bigr]
	< e^{-t^2 / 4}$$
for all $t > 0$. This is comparable to the result of Theorem
\ref{eigenvalues} for $\lambda_1(X)$ and $\lambda_n(X)$, but the
same level of concentration holds for eigenvalues in the bulk of the spectrum,
which is not the case in Theorem \ref{eigenvalues}.

The result of Theorem \ref{eigenvalues} for $\lambda_1(X)$ and
$\lambda_n(X)$ (stated in less generality) was shown by
Krivelevich and Vu in \cite{KV}.
After a preliminary version of this paper
was written, we learned that Alon, Krivelevich, and Vu \cite{AKV} showed
that for $1 \le k \le n$,
$$ \P \bigl[ |\lambda_k(X) - \M \lambda_k(X) | \ge t \bigr] \le
    4 \exp\left[ -\frac{t^2}{8 k^2 D^2} \right]$$
for all $t > 0$, and that the same holds if $\lambda_k(X)$ is replaced by
$\lambda_{n-k+1}(X)$. The approach in \cite{AKV} handles the lack of
convexity of $\lambda_k$ by not using the $q = 2$ case of Corollary
\ref{l_p-con}, but instead applying Theorem \ref{Talagrand} by directly
estimating the convex hull
distances involved. Our Theorem \ref{eigenvalues} improves the order of
fluctuations of $\lambda_k(X)$ from $O(k)$ (as in \cite{AKV}) to
$O(\sqrt{k})$. It is also conjectured in \cite{AKV} that $\lambda_k(X)$ 
should be concentrated at least as strongly as $\lambda_1(X)$, as one
obtains from Theorem \ref{Gaussian} in the Gaussian case. We emphasize
again that we are dealing only with large deviations here. As we have
already indicated in the introduction, the tail estimate 
\eqref{lambda_1-ineq} for the extreme eigenvalues is not sharp for
$t = o(\sqrt{n})$; furthermore, it is likely that concentration is
even tighter for eigenvalues in the bulk of the spectrum.

It follows as in the discussion following Corollary \ref{l_p-con}
that Theorem \ref{eigenvalues} implies that $\E \lambda_k(X)$ differs
by at most $O(\sqrt{k})$ from the number $M_k$ which appears in the
statement of the theorem. One can also show that the number $M_k$
which appears in the statement of the theorem differs by at most
$O(\sqrt{k})$ from $\M \lambda_k(X)$. By using the separate bounds for
deviations above and below the median in the situation of Corollary
\ref{l_p-con}, we have
$$|M_k - \M \lambda_k(X) | \le 2\sqrt{6 \log 2}(\sqrt{k}
+ \sqrt{k-1}) D.$$

We can also obtain a similar result to Theorem \ref{eigenvalues} for
singular values in the rectangular case. Let $l = \min \{ m, n \}$.
For an $m \times n$ matrix $A$, we denote by $s_1(A) \ge s_2(A) \ge \cdots
\ge s_l(A) \ge 0$ the singular values of $A$, counted with multiplicity;
that is, $s_k(A) = \lambda_k \bigl( (A^* A)^{1/2} \bigr)$.

\begin{thm} \label{singular}
Suppose the entries $x_{jk}$ of $X$ are independent complex
random variables, each satisfying the condition (\ref{r-i}) in the
statement of Theorem \ref{eigenvalues}. Then
$$ \P \bigl[ | s_1(X) - \M s_1(X) | \ge t \bigr] \le 4 e^{-t^2 / 4 D^2}$$
for all $t > 0$. Furthermore, for each $2 \le k \le 
\min \{ m, n \}$, there exists an $M_k \in \R$ such that
$$\P \bigl[ | s_k(X) - M_k | \ge t \bigr] \le
    8 \exp \left[ - \left(\frac{t}{2(\sqrt{k} + \sqrt{k-1}) D }\right)^2
    \right] \le 8 \exp\left[-\frac{t^2}{16 k D^2}\right]$$
for all $t > 0$.
\end{thm}

The proof is similar to the proof of Theorem \ref{eigenvalues}, using
in place of the functions $F_k$ the Ky Fan $k$-norms, defined by
$$\| A \|_{(k)} = \sum_{j=1}^k s_j(A)$$
for $1 \le k \le \min \{ m, n \}$. We remark that the triangle inequality,
and hence convexity, for the Ky Fan norms can be proved as a formal
consequence of the convexity of the functions $F_k$.


\section*{Acknowledgments}
This paper is part of the author's Ph.D. thesis, written under the
supervision of Profs. S. Szarek and E. Werner. The author wishes to
thank Prof. S. Szarek for many valuable discussions.


\end{document}